\documentclass[a4paper,11pt, reqno]{amsart}
\usepackage{amssymb,amsthm,amsmath}
 \usepackage{enumerate}
\usepackage{ifthen}
\usepackage{graphicx}
\usepackage{cite}
\usepackage{tikz} 
\usetikzlibrary{arrows.meta}
\usepackage{amsmath,amscd,amssymb}
\usepackage{latexsym}
\usepackage[colorlinks,citecolor=blue,pagebackref,hypertexnames=false]{hyperref}

\numberwithin{equation}{section}
\allowdisplaybreaks[2]
\theoremstyle{plain}
\newtheorem{theorem}{Theorem}[section]

\newtheorem{proposition}[theorem]{Proposition}
\newtheorem{conjecture}[theorem]{Conjecture}

\theoremstyle{definition}

\theoremstyle{remark}
\newtheorem{remark}[theorem]{Remark}

\newtheorem{notation}[theorem]{Notation}

\newtheorem{case[theorem]}{Case}

\def\norm#1.#2.{\lVert#1\rVert_{#2}}

\title[Strichartz estimates at the critical summability exponent]{Strichartz estimates involving orthonormal systems at the critical summability exponent
}
\author{Guoxia Feng}
\author{Manli Song}
\author{Huoxiong Wu}
\address{\endgraf School of Mathematical Sciences, Xiamen University, Xiamen 361005, China}
\email{gxfeng@mail.nwpu.edu.cn}

\address{\endgraf School of Mathematics and Statistics, Northwestern Polytechnical University, Xi'an, Shaanxi 710129, China}
\email{mlsong@nwpu.edu.cn}

\address{\endgraf School of Mathematical Sciences, Xiamen University, Xiamen 361005, China}
\email{huoxwu@xmu.edu.cn}

\keywords{Strichartz estimate for orthonormal functions, Schr\"odinger equation}
\subjclass[2010]{35Q41 (primary), 42B37, 35B65 (secondary)}
\date{today}

\begin{document}
\maketitle

	\allowdisplaybreaks

\begin{abstract}
The primary objective of this paper is to investigate the orthonormal Strichartz estimates at the critical summability exponent for the Schrödinger operator $e^{it\Delta}$ with initial data from the homogeneous Sobolev space $\dot{H}^s (\mathbb{R}^n)$. We prove new global strong-type orthonormal Strichartz estimates in the interior of $ODCA$ at the optimal summability exponent $\alpha=q$, thereby substantially supplymenting the work of Bez-Hong-Lee-Nakamura-Sawano \cite{Bez-Hong-Lee-Nakamura-Sawano}. Our approach is based on restricted weak-type orthonormal estimates, real interpolation argument and the advantageous condition $q<p$ in the interior of $ODCA$.
\end{abstract}
	\tableofcontents 	

\section{Introduction}
\subsection{Single-function Strichartz estimates}
The frequency-localised Strichartz estimates involving a single function address the following space-time estimates
\begin{equation}\label{local}
\|e^{it\Delta}P_0f\|_{L^{2q}(\mathbb{R},L^{2p}(\mathbb{R}^n))}\lesssim \|f\|_{L^2(\mathbb{R}^n)}
\end{equation}
for all spatial dimensions $n\geq1$ and the admissible pairs $(p,q)$ satisfying the condition \begin{equation}\label{AC}
1\leq p,q\leq \infty,\; \frac{2}{q}+\frac{n}{p}\leq n, \text{ and } (q,p,n)\neq (1,\infty,2),
\end{equation}
where $e^{it\Delta}$, $P_k$, $k\in\mathbb{Z}$ and $L^{2q}(\mathbb{R},L^{2p}(\mathbb{R}^n))$ are defined in Section \ref{Preliminaries}. Here, when $\frac{2}{q}+\frac{n}{p}\leq n$ holds, we shall refer to it as the sharp admissible case, when $\frac{2}{q}+\frac{n}{p}\leq n$ as the non-sharp admissible case, and when $q=1,\infty$ or $p=\infty$ as boundary cases.

The estimates \eqref{local} for case $q=p$ were initially proved in the pioneering work by Strichartz \cite{Str} and then the other cases were deduced from a combination of an abstract functional analysis argument known as the $TT^*$ duality argument and an $L^1\rightarrow L^\infty$ dispersive estimate, except the endpoint case $(q,p)=(1,\frac{n}{n-2})$ for $n\geq3$, which was proved by Keel-Tao \cite{KT}. 

By a scaling argument, one can get from \eqref{local} that
\begin{equation*}
\|e^{it\Delta}(-\Delta)^{-\frac{s}{2}}P_kf\|_{L^{2q}(\mathbb{R},L^{2p}(\mathbb{R}^n))}\lesssim 2^{-\frac{k}{2}(n-2s-\frac{2}{q}-\frac{n}{p})}\|f\|_{L^2(\mathbb{R}^n)}.
\end{equation*}
Then we can upgrade the frequency-localised estimates to frequency-global cases by putting together these estimates for each dyadic piece by the Minkowski's inequality and the Littlewood-Paley square function theorem, 
\begin{equation}\label{global}
\|e^{it\Delta}f\|_{L^{2q}(\mathbb{R},L^{2p}(\mathbb{R}^n))}\lesssim \|f\|_{\dot{H}^s(\mathbb{R}^n)},
\end{equation}
if $(p,q)$ satisfies the admissible condition \eqref{AC}, $p\neq \infty$ and the necessary scaling condition
\begin{equation*}
\frac{2}{q}+\frac{n}{p}=n-2s,
\end{equation*}
and we can see that $s\in [0,\frac{n}{2})$.
On the other hand, for the boundary case $p=\infty$, the estimate \eqref{global}, namely,
\begin{equation}\label{boundary}
\|e^{it\Delta}f\|_{L^{2q}(\mathbb{R},L^{\infty}(\mathbb{R}^n))}\lesssim \|f\|_{\dot{H}^{\frac{n}{2}-\frac{1}{q}}(\mathbb{R}^n)}
\end{equation}
 needs special treatment due to the failure of the Littlewood-Paley square function theorem in $L^\infty$. It was very recently completed solved by Guo-Li-Nakanishi-Yan \cite{Guo-Li-Nakanishi-Yan}. In particular, it is shown that the estimate in \eqref{boundary} holds for $1<q<\infty$ when $n\geq2$ and holds for $2\leq q<\infty$ when $n=1$. But it fails for $q=\infty$, which follows from the well-known failure of the corresponding Sobolev embedding estimate, and also fails for $q=1$ when $n\geq2$ (when $n\geq3$, see \cite{Guo-Li-Nakanishi-Yan} and when $n=2$, see \cite{Montgomery-Smith}). 
\subsection{Strichartz estimates involving orthonormal systems}
 In recent years, significant attention has been devoted by numerous researchers to the extension from single-function Strichartz estimates  to versions involving systems of orthonormal functions taking the form
    \begin{equation}\label{abstract}
\left\|\sum_{j\in J}\lambda_j|e^{it L}f_j|^2\right\|_{L^q(I,L^p(X))} \lesssim \|\{\lambda_j\}_{j\in J}\|_{\ell^{\alpha}},
	\end{equation}
	for families of orthonormal data  $\{f_j\}_{j=1}^\infty$ in a given Hilbert space, such as the Lebesgue space $L^2(X)$ or the homogeneous Sobolev space $\dot{H}^s(X)$. Here, $I\subseteq\mathbb{R}$ is an interval, $L$ is a non-negative self-adjoint operator on $L^2(X)$, where $(X,d\mu)$ is a given measure space, $\{\lambda_j\}_{j=1}^\infty \in \ell^{\alpha}$  and $ q, p, \alpha \in [1, \infty]$ satisfying some particular conditions. 
    
The estimates \eqref{abstract} first appeared in the pioneering work of Frank-Lewin-Lieb-Seiringer \cite{FLLS}, which generalized single-function Strichartz estimates \eqref{global} of the Schr\"odinger operator $e^{it\Delta}$ to systems of orthonormal functions. Later, Frank-Sabin \cite{FS, FS1}, Bez-Hong-Lee-Nakamura-Sawano \cite{Bez-Hong-Lee-Nakamura-Sawano} and Bez-Lee-Nakamura \cite{BEZ} investigated a wide class of dispersive equations for the Laplacian on the Euclidean space and established orthonormal Strichartz estimates in the case of wave, Klein-Gordon and fractional Schr\"odinger operators.  This also works for dispersive equations with potentials, see \cite{Hoshiya2,Hoshiya, Feng-Song, FMSW,GMS,Nguyen}. See \cite{MS,shyam, fmsw} for the orthonormal Strichartz estimates for the Dunkl Laplacian or (k,a)-generalized Laguerre operator and \cite{nakamura, FS2, Ren-Zhang, Wang-Zhang-Zhang} on compact Riemannian manfiolds.

In our article, we focus on the orthonormal Strichartz estimates of the form \eqref{abstract} for the classical Schr\"odinger operator $e^{it\Delta}$. We now briefly summarize some particular known results.   In order to state precisely, first we establish some notation.
\begin{notation}
We introduce the following points:
\begin{equation*}
\begin{array}{llllll}
&O = (0,0),\quad & A = (\frac{n-1}{n+1},\frac{n}{n+1}), \quad & B = (1,0), \quad & C = (0,1), \\
& D = (\frac{n-2}{n},1), & E = (\frac{n-1}{2n},\frac{1}{2}), \quad & F  = (\frac{n}{n+2},\frac{n}{n+2}).
\end{array}
\end{equation*}  

For $n\geq 3$, see Figure 1. The line segment $[B,D]$ corresponds to the sharp admissile case. The interior of $OBDC$ corresponds the non-sharp admissile case. The line segment $[O,C]$, $[C,D)$ and $[O,B)$ correspond to boundary cases.
\begin{center}\begin{tikzpicture} 
\draw[thick,->, >=latex] (0,0) -- (0,7) node[left] {\(\frac{1}{q}\)};\draw[thick,->, >=latex] (-5/2,0) -- (7,0) node[below] {\(\frac{1}{p}\)};\node at (0,6) [above right]{\(1\)};\node at (6,0) [below]{\(1\)};\node at (0,0) [below]{\(O\)};\node at (0,6) [above left]{\(C\)};\node at (2,6) [above]{\(D\)};\node at (6,0) [above right]{\(B\)};\node at (18/5,18/5) [above right]{\(F\)}; \node at (3,9/2) [above right]{\(A\)};\node at (2,3) [above left]{\(E\)};\node at (0,3) [left]{\(\frac{1}{2}\)};\node at (2,0) [below]{\(\frac{n-1}{2n}\)};\node at (0,9/2) [left]{\(\frac{n}{n+1}\)};\node at (3,0) [below]{\(\frac{n-1}{n+1}\)};\node at (0,18/5) [left]{\(\frac{n}{n+2}\)};\node at (18/5,0) [below]{\(\frac{n}{n+2}\)};\node at (3,-1) [below] {\(\text{Figure 1: }n\geq3\)};\draw[thick] (6,0) -- (6,6);\draw[thick] (6,0) -- (2,6);\draw[thick] (0,6) -- (6,6);\draw[thick] (0,0) -- (3,9/2);\draw[dotted] (0,3)--(2,3);\draw[dotted] (2,0)--(2,3);\draw[dotted] (3,0)--(3,9/2);\draw[dotted] (0,9/2)--(3,9/2);\draw[dotted] (18/5,0)--(18/5,18/5);\draw[dotted] (0,18/5)--(18/5,18/5);\end{tikzpicture}\end{center}

For $n=2$, see Figure 2. The points $C$ and $D$ coincide. In this case, the quadrilateral $OCDA$ reduces to the triangle $OAD$ and the line segment $[O,C]$ is exactly $[O,D]$.
\begin{center}\begin{tikzpicture}\draw[thick,->, >=latex] (0,0) -- (0,7) node[left] {\(\frac{1}{q}\)};\draw[thick,->, >=latex] (-5/2,0) -- (7,0) node[below] {\(\frac{1}{p}\)};\node at (0,6) [above right]{\(1\)};\node at (6,0) [below]{\(1\)};\node at (0,0) [below]{\(O\)};\node at (0,6) [above left]{\(D\)};\node at (6,0) [above right]{\(B\)};\node at (3,3) [above right]{\(F\)};\node at (2,4) [above right]{\(A\)};\node at (3/2,3) [above left]{\(E\)};\node at (0,3) [left]{\(\frac{1}{2}\)};\node at (3/2,0) [below]{\(\frac{1}{4}\)};\node at (0,4) [left]{\(\frac{2}{3}\)};\node at (2,0) [below]{\(\frac{1}{3}\)};\node at (0,3) [left]{\(\frac{1}{2}\)};\node at (3,0) [below]{\(\frac{1}{2}\)};\node at (3,-1) [below] {\(\text{Figure 2: }n=2\)};\draw[thick] (6,0) -- (6,6);\draw[thick] (6,0) -- (0,6);\draw[thick] (0,6) -- (6,6);\draw[thick] (0,0) -- (2,4);\draw[dotted] (0,3)--(3/2,3);\draw[dotted] (3/2,0)--(3/2,3);\draw[dotted] (2,0)--(2,4);\draw[dotted] (0,4)--(2,4);\draw[dotted] (3,0)--(3,3);\draw[dotted] (0,3)--(3,3);\end{tikzpicture}\end{center}

For $n=1$, see Figure 3. The points $A$ and $E$ coincide. 
\begin{center}
\begin{tikzpicture}
            \draw[thick,->, >=latex] (0,0) -- (0,7) node[left] {\(\frac{1}{q}\)};
            \draw[thick,->, >=latex] (-5/2,0) -- (7,0) node[below] {\(\frac{1}{p}\)};
            \node at (0,6) [above right]{\(1\)};
            \node at (6,0) [below]{\(1\)};
            \node at (0,0) [below]{\(O\)};
            \node at (0,6) [above left]{\(C\)};
            \node at (-6,6) [above]{\(D\)};
            \node at (6,0) [above right]{\(B\)};
            \node at (2,2) [above right]{\(F\)};
            \node at (0,3) [above right]{\(A\)};
            \node at (0,2) [left]{\(\frac{1}{3}\)};
            \node at (0,3) [below left]{\(\frac{1}{2}\)};
            \node at (2,0) [below]{\(\frac{1}{3}\)};
            \node at (3,-1) [below] {\(\text{Figure 3: }n=1\)};
            \draw[thick] (0,0) -- (2,2);
            \draw[thick] (6,0) -- (6,6);
            \draw[thick] (6,0) -- (0,3);
            \draw[dotted] (-6,6) -- (0,3);
            \draw[thick] (0,6) -- (6,6);
            \draw[dotted] (0,6) -- (-6,6);
            \draw[dotted] (2,0) -- (2,2);
            \draw[dotted] (0,2)--(2,2);
        \end{tikzpicture}
\end{center}

For points $X_j \in \mathbb{R}^2$, $j=1,2,3,4$, we write
\begin{align*}
[X_1,X_2] & = \{ (1-t)X_1 + tX_2 : t \in [0,1]\} \\
[X_1,X_2) & = \{ (1-t)X_1 + tX_2 : t \in [0,1)\} \\
(X_1,X_2) & = \{ (1-t)X_1 + tX_2 : t \in (0,1)\}
\end{align*}
for line segments connecting $X_1$ and $X_2$, including or excluding $X_1$ and $X_2$ as appropriate. We write
$
X_1X_2X_3
$
for the convex hull of $X_1, X_2, X_3$, and
$
{\rm int}\, X_1X_2X_3
$
for the interior of $X_1X_2X_3$. Similarly,
$
X_1X_2X_3X_4
$
denotes the convex hull of $X_1, X_2, X_3,X_4$, and
$
{\rm int}\, X_1X_2X_3X_4
$
denotes the interior of $X_1X_2X_3X_4$. In particular,
\[
{\rm int}\, OAB = \bigg\{ \bigg(\frac{1}{p},\frac{1}{q}\bigg) \in (0,1)^2 : \text{$\frac{1}{q} < \frac{n}{(n-1)p}$ and $\frac{2}{q} + \frac{n}{p} < n$}\bigg\}
\]
and
\[
{\rm int}\, OCDA = \bigg\{ \bigg(\frac{1}{p},\frac{1}{q}\bigg) \in (0,1)^2 : \text{$\frac{1}{q} >\frac{n}{(n-1)p}$ and $\frac{2}{q} + \frac{n}{p} < n$}\bigg\}.
\]

%We remark that the line segment $[B,A)$ corresponds to the range of estimates in Theorem \ref{M-S}.

Throughout the paper, we need the exponent $\alpha^*(p,q)$ determined by
\[
\frac{n}{\alpha^*(p,q)}=\frac{1}{q}+\frac{n}{p}.
\]
Note that if $\frac{2}{q} + \frac{n}{p} = n$ (corresponding to the case $s=0$) then $\alpha^*(p,q) = \frac{2p}{p+1}$. Also, note that if $(\frac{1}{p},\frac{1}{q})$ belongs to the line segment $[O,A]$, then $\alpha^*(p,q) = q$, and if $(\frac{1}{p},\frac{1}{q})$ belongs to the line segment $[O,B]$, then $\alpha^*(p,q) = p$. In addition, if $(\frac{1}{p},\frac{1}{q})$ belongs to the region int $OCDA$, then $\alpha^*(p,q)<q<p$. If $(\frac{1}{p},\frac{1}{q})$ belongs to the region int $OAB$, then $\alpha^*(p,q)>q$.
\end{notation}
Now we restate one-single frequency-localised Strichartz estimates \eqref{local}.
\begin{itemize}
\item For $n\geq 3$, the estimate \eqref{local} holds for $(\frac{1}{p},\frac{1}{q})$ belonging to $OBDC$.
\item For $n=2$, the estimate \eqref{local}holds for $(\frac{1}{p},\frac{1}{q})$ belonging to $OBD\setminus D$, but fails at points $D$.
\item For $n=1$, the estimate \eqref{local} holds for $(\frac{1}{p},\frac{1}{q})$ belonging to $OAB$.
\end{itemize}

We also need to restate one-single frequency-global Strichartz estimates \eqref{global} with the necessary scaling condition $\frac{2}{q}+\frac{n}{p}=n-2s$.
\begin{itemize}
\item For $n\geq 2$, the estimate \eqref{global} holds for $(\frac{1}{p},\frac{1}{q})$ belonging to $OBDC\setminus \{O,C\}$, but fails at points $O$ and $C$.
\item For $n=1$, the estimate \eqref{global} holds for $(\frac{1}{p},\frac{1}{q})$ belonging to $OAB\setminus O$, but fails at point $O$.
\end{itemize}

Recently, single-function Strichartz estimates have been substantially generalized to the context of orthonormal systems in $L^2(\mathbb{R}^n)$ by Frank-Lewin-Lieb-Seiringer \cite{FLLS} and Frank-Sabin \cite{FS}.
\begin{theorem}[\cite{FLLS,FS}]\label{Lebesgue} Suppose $n\geq1$. If $(\frac{1}{p},\frac{1}{q})$ belongs to $(A,B]$ and $\alpha=\frac{2p}{p+1}$, then
\begin{equation}\label{lebesgue}
			\bigg\|\sum_{j\in J} \lambda_j|e^{it\Delta}f_j|^2\bigg\|_{L^q(\mathbb{R}, L^p(\mathbb{R}^n))}\lesssim \|\{\lambda_j\}_{j\in J}\|_{\ell^{\alpha}},
\end{equation}	
  holds for all families of orthonormal functions $\{f_j\}_{j\in J}$ in $L^2(\mathbb{R}^n)$ and all sequences $\{\lambda_j\}_{j\in J}$ in $\ell^{\alpha}$. This is sharp in the sense that the estimate \eqref{lebesgue} fails for all $\alpha>\frac{2p}{p+1}$. Furthermore, for $n\geq2$, the estimate \eqref{lebesgue} holds at the endpoint $A$ for all $\alpha<\frac{2p}{p+1}$ and fails for $\alpha=\frac{2p}{p+1}$.
\end{theorem}
It is naturally to conjecture that for $n=1$ the estimate \eqref{lebesgue} should hold at the endpoint $A$ for  $\alpha<2$. However, as far as we are aware, whether it holds or not is a challenging open problem. Bez-Lee-Nakamura \cite{Bez-Lee-Nakamura1} obtained the optimal range of $\alpha$ for a weak-type version of \eqref{lebesgue}, while it seems to be non-trivial to upgrade this to a strong-type estimate. 
\begin{theorem}[\cite{Bez-Lee-Nakamura1}]
The estimate 
\begin{equation*}
			\bigg\|\sum_{j\in J} \lambda_j|e^{it\partial_x^2}f_j|^2\bigg\|_{L^{2,\infty}(\mathbb{R}, L^\infty(\mathbb{R}))}\lesssim \|\{\lambda_j\}_{j\in J}\|_{\ell^{\alpha}},
\end{equation*}	
holds for all systems of orthonormal functions $\{f_j\}_{j\in J}$ in $L^2(\mathbb{R})$ and all sequences $\{\lambda_j\}_{j\in J}$ in $\ell^{\alpha}$ if and only if $\alpha<2$.
\end{theorem}
\begin{remark}
The strong-type estimate 
\begin{equation*}
			\bigg\|\sum_{j\in J} \lambda_j|e^{it\partial_x^2}f_j|^2\bigg\|_{L^{2}(\mathbb{R}, L^\infty(\mathbb{R}))}\lesssim \|\{\lambda_j\}_{j\in J}\|_{\ell^{\alpha}},
\end{equation*}	
was observed by Bez-Lee-Nakamura \cite{Bez-Lee-Nakamura1} to hold in a smaller range $\alpha\leq \frac{4}{3}.$
\end{remark}
Returning to Theorem \ref{Lebesgue}, for $n\geq3$, note that it follows immediately from the triangle inequality and one-single function Strichartz estimate \eqref{global} for the sharp-admissible case $(\frac{1}{p},\frac{1}{q})\in [B,D]$ that
\begin{equation}\label{trivial}
			\bigg\|\sum_{j\in J} \lambda_j|e^{it\Delta}f_j|^2\bigg\|_{L^q(\mathbb{R}, L^p(\mathbb{R}^n))}\leq \sum_{j\in J} |\lambda_j|\big\|e^{it\Delta}f_j|^2\big\|_{L^q(\mathbb{R}, L^p(\mathbb{R}^n))}\lesssim \|\{\lambda_j\}_{j\in J}\|_{\ell^1},
\end{equation}	
which gives \eqref{lebesgue} with $\alpha=1$ without making use of the orthogonality and the pertinent point here is to raise $\alpha$ as far as possible by capitalising on the orthogonality of $\{f_j\}_{j\in J}$. At the Keel-Tao endpoint $D$, the exponent $\alpha=1$ can not be improved (see \cite{FS1}). We see that the point $A$ plays a crucial role in extending the form of estimates \eqref{lebesgue} on the line segment $[B,D]$. Indeed, interpolating \eqref{lebesgue} between the point $(\frac{1}{p},\frac{1}{q})$ on the line $(A,B]$ arbitrarily close to $A$ (with $\alpha=\frac{2p}{p+1}$) and the Keel-Tao endpoint $D$ (with $\alpha=1$), it implies the estimate \eqref{lebesgue} for all points $(\frac{1}{p},\frac{1}{q})\in [A,D)$ with $\alpha<q$. Also a remarkable phenomena here is that the sharp value of $\alpha$ reaches its maximum at the point $A$. Therefore, $A$ may be considered as an endpoint case. For $n=2$, interpolating \eqref{lebesgue} between the point $(\frac{1}{p},\frac{1}{q})\in (A,B]$ arbitrarily close to $A$ (with $\alpha=\frac{2p}{p+1}$) and the point $(\frac{1}{p},\frac{1}{q})\in [A,D)$ arbitrarily close to the Keel-Tao endpoint $D$ (with $\alpha=1$), it implies the estimate \eqref{lebesgue} for all points $(\frac{1}{p},\frac{1}{q})\in [A,D)$ with $\alpha<q$. On the other hand, it was proved in \cite{FS1} that, for points $(\frac{1}{p},\frac{1}{q})\in [A,D)$, the estimate \eqref{lebesgue} fails when $\alpha>q$. 
\begin{theorem}[\cite{FS1}]
Suppose $n\geq2$. If $(\frac{1}{p},\frac{1}{q})$ belongs to $[A,D)$, then for any $\alpha<q$
\begin{equation}\label{interpolation}
			\bigg\|\sum_{j\in J} \lambda_j|e^{it\Delta}f_j|^2\bigg\|_{L^q(\mathbb{R}, L^p(\mathbb{R}^n))}\lesssim \|\{\lambda_j\}_{j\in J}\|_{\ell^{\alpha}},
\end{equation}	
  holds for all families of orthonormal functions $\{f_j\}_{j\in J}$ in $L^2(\mathbb{R}^n)$ and all sequences $\{\lambda_j\}_{j\in J}$ in $\ell^{\alpha}$. This is sharp in the sense that the estimate \eqref{lebesgue} fails for all $\alpha>q$. 
\end{theorem}
Therefore, for $n\geq2$, the only remaining issue is the critical case $\alpha=q$. The following interesting conjecture was raised in \cite{FLLS} (see also in \cite{FS1}).
\begin{conjecture}[\cite{FLLS,FS1}]
 Let $n\geq1$. At the endpoint $A=(\frac{1}{p},\frac{1}{q})=(\frac{n-1}{n+1},\frac{n}{n+1})$, the restricted-type estimate
 \begin{equation}\label{conjecture}
			\bigg\|\sum_{j\in J} \lambda_j|e^{it\Delta}f_j|^2\bigg\|_{L^q(\mathbb{R}, L^p(\mathbb{R}^n))}\lesssim \|\{\lambda_j\}_{j\in J}\|_{\ell^{q,1}},
\end{equation}
 holds for all families of orthonormal functions $\{f_j\}_{j\in J}$ in $L^2(\mathbb{R}^n)$ and all sequences $\{\lambda_j\}_{j\in J}$ in $\ell^{q,1}$.
\end{conjecture}
For $n\geq3$, if the conjectured estimates \eqref{conjecture} were true, proceeding via real interpolation with the estimate \eqref{trivial} at the Keel-Tao endpoint, it gives the critical estimates \eqref{interpolation} for $\alpha=q$. Interestingly, such estimates fail for $n=1$. In fact, this was proved by Bez-Hong-Lee-Nakamura-Sawano \cite{Bez-Hong-Lee-Nakamura-Sawano} using a semi-classical limit argument and demonstrating the failure of certain induced estimates through a geometric argument that relied on the existence of Nikodym sets with zero Lebesgue measure (see Theorem $3,5$ in \cite{Bez-Hong-Lee-Nakamura-Sawano}). It remains open whether the conjectured estimates \eqref{conjecture} hold for $n\geq2$.

At the same time, Bez-Hong-Lee-Nakamura-Sawano \cite{Bez-Hong-Lee-Nakamura-Sawano} also considered the case of orthonormal system in the homogeneous Sobolev space $\dot{H}^s(\mathbb{R}^n)$ (which corresponds to the non-sharp adimissble case). Firstly, they established the following optimal frequency-localised estimates. 
\begin{theorem}[\cite{Bez-Hong-Lee-Nakamura-Sawano}] \label{frequency-localized} Suppose $n\geq3$.\\

$(1)$ If $(\frac{1}{p},\frac{1}{q})$ belongs to $ OAB\setminus A$  and $\alpha=\alpha^*(p,q)$, then
\begin{equation}\label{OAB}
			\bigg\|\sum_{j\in J}\lambda_j|e^{it\Delta}P_0f_j|^2\bigg\|_{L^q(\mathbb{R}, L^p(\mathbb{R}^n))}\lesssim_\psi \|\{\lambda_j\}_{j\in J}\|_{\ell^{\alpha}},
\end{equation}	
  holds for all families of orthonormal functions $\{f_j\}_{j\in J}$ in $L^2(\mathbb{R}^n)$ and all sequences $\{\lambda_j\}_{j\in J}$ in $\ell^{\alpha}$.  This is sharp in the sense that the estimate fails for all $\alpha>\alpha^*(p,q)$.
  
$(2)$ If $(\frac{1}{p},\frac{1}{q})$ belongs to $OCDA\setminus [A,D)$ and $\alpha=q$, then the results \eqref{OAB} hold true. This is sharp in the sense that the estimate fails for all $\alpha>q$.
\end{theorem}
\begin{remark}
Indeed, for $n\geq 3$, Bez-Hong-Lee-Nakamura-Sawano \cite{Bez-Hong-Lee-Nakamura-Sawano} first established the frequency-localised estimates for four special cases: the origin point $O$, the line segments $[C,D]$, $(A,E]$ and $[A,B]$. All other cases then followed via complex interpolation arguments.
\end{remark}

For $1\leq n\leq 2$, the corresponding frequency-localized estimates are provided in the following remarks.
\begin{itemize}
    \item \label{n=2} Let $n=2$.  Then the points $C$ and $D$ coincide. If $(\frac{1}{p},\frac{1}{q})$ belongs to $ OAB\setminus A$  and $\alpha=\alpha^*(p,q)$, the results \eqref{OAB} hold true. Moreover,  if $(\frac{1}{p},\frac{1}{q})$ belongs to $OAD\setminus [A,D]$ with $\alpha=q$, the results \eqref{OAB} still hold.
    \item \label{n=1} Let $n=1$, then the points $A$ and $E$ coincide. If  $(\frac{1}{p},\frac{1}{q})$ belongs to $ OAB\setminus (O,A]$ and $\alpha=\alpha^*(p,q)$, then the results in \eqref{OAB} holds true.
\end{itemize}

For the general case $k\in\mathbb{Z}$, by a rescaling argument, it follows immediately from the estimate \eqref{OAB} in Theorem \ref{frequency-localized} that
\begin{equation}\label{general-k}
			\bigg\|\sum_{j\in J}\lambda_j|e^{it\Delta}(-\Delta)^{-\frac{s}{2}}P_kf_j|^2\bigg\|_{L^q(\mathbb{R}, L^p(\mathbb{R}^n))}\lesssim_\psi 2^{k(n-2s-\frac{2}{q}-\frac{n}{p})}\|\{\lambda_j\}_{j\in J}\|_{\ell^{\alpha}},
\end{equation}
holds for all families of orthonormal functions $\{f_j\}_{j\in J}$ in $L^2(\mathbb{R}^n)$ and all sequences $\{\lambda_j\}_{j\in J}$ in $\ell^{\alpha}$. 

Unfortunately, the frequency-localized estimates in Theorem \ref{frequency-localized} do not readily extend to general data, unlike their counterparts in the single-function case. By Proposition \ref{upgrading} in Section \ref{LP}, we can sum up the frequency-localised estimates \eqref{general-k} and upgrade them to restricted weak-type frequency-global results. The following strong-type global frequency estimates can be derived through a series of successive real and complex interpolation arguments.
\begin{theorem}[\cite{Bez-Hong-Lee-Nakamura-Sawano}]\label{globally}
$(1)$ Let $n\geq 1$. If $(\frac{1}{p},\frac{1}{q})$ belongs to int $OAB$, $\alpha=\alpha^*(p,q)$ and $2s=n-\left(\frac{2}{q}+\frac{n}{p}\right)$, then
\begin{equation}\label{global-orthonormal}
			\bigg\|\sum_{j\in J} \lambda_j|e^{it\Delta}f_j|^2\bigg\|_{L^q(\mathbb{R}, L^p(\mathbb{R}^n))}\lesssim \|\{\lambda_j\}_{j\in J}\|_{\ell^{\alpha}}, 
\end{equation}	
  holds for all families of orthonormal functions $\{f_j\}_{j\in J}$ in $\dot{H}^{s}(\mathbb{R}^n)$ and all sequences $\{\lambda_j\}_{j\in J}$ in $\ell^{\alpha}$. This is sharp in the sense that the estimate fails for all $\alpha>\alpha^*(p,q)$.
  
$(2)$ Let $n\geq 2$. If $(\frac{1}{p},\frac{1}{q})$ belongs to int $OCDA$, $2s=n-\left(\frac{2}{q}+\frac{n}{p}\right)$ and $\alpha<q$, then the estimate \eqref{global-orthonormal} holds for all families of orthonormal functions $\{f_j\}_{j\in J}$ in $\dot{H}^{s}(\mathbb{R}^n)$ and all sequences $\{\lambda_j\}_{j\in J}$ in $\ell^{\alpha}$. This is sharp in the sense that the estimate fails for all $\alpha>q$. 
\end{theorem}
It was also shown in \cite{Bez-Hong-Lee-Nakamura-Sawano} that 
\begin{equation*}
    \alpha\leq \min\{q, \alpha^*(p,q)\}
\end{equation*}
is a necessary condition for the estimate \eqref{global-orthonormal} to hold and hence Theorem \ref{globally} is close to optimal except for the critical summability exponent $\alpha=q$. We shall present our new results in this direction below in Section \ref{New}.

Next we discuss the boundary cases. For $n\geq 3$, as \eqref{trivial}, the orthonormal Strichartz estimates at the critical summability exponent $\alpha=1$ on the boundary line $[C,D)$ follows immediately from the triangle inequality and the one-single function Strichartz estimate \eqref{global}. For the boundary case $q=\infty$ (corresponds to the line segment $(O,B)$), $\alpha^*(p,\infty)=p$, Bez-Hong-Lee-Nakamura-Sawano \cite{Bez-Hong-Lee-Nakamura-Sawano} obtained the following restricted weak-type estimate.
\begin{proposition}[\cite{Bez-Hong-Lee-Nakamura-Sawano}]
    Let $n\geq1$. If $1<p<\infty$ and $2s=n-\frac{n}{p}$, then the estimate
    \begin{equation}\label{orthonormal-Boundary}
			\bigg\|\sum_{j\in J} \lambda_j|e^{it\Delta}f_j|^2\bigg\|_{L^\infty(\mathbb{R}, L^{p,\infty}(\mathbb{R}^n))}\lesssim \|\{\lambda_j\}_{j\in J}\|_{\ell^{p,1}}, 
\end{equation}	
holds for all families of orthonormal functions $\{f_j\}_{j\in J}$ in $\dot{H}^{s}(\mathbb{R}^n)$ and all sequences $\{\lambda_j\}_{j\in J}$ in $\ell^{p,1}$. 
\end{proposition}
It is an interesting open problem whether it is possible to upgrade the weak-$L^p$ norm at the left hand to $L^p$.

For the boundary case $p=\infty$ (corresponds to the line segment $(O,C)$), Bez-Kinoshita-Shiraki \cite{Bez-Kinoshita-Shiraki} extended \eqref{boundary} to the estimate involving orthonormal systems.

\begin{theorem}[\cite{Bez-Kinoshita-Shiraki}] The estimate 
   \begin{equation*}
			\bigg\|\sum_{j\in J} \lambda_j|e^{it\Delta}f_j|^2\bigg\|_{L^q(\mathbb{R}, L^\infty(\mathbb{R}^n))}\lesssim \|\{\lambda_j\}_{j\in J}\|_{\ell^{\alpha}}, 
\end{equation*}	
holds for all families of orthonormal functions $\{f_j\}_{j\in J}$ in $\dot{H}^{s}(\mathbb{R}^n)$ with $s=\frac{n}{2}-\frac{1}{q}$ and all sequences $\{\lambda_j\}_{j\in J}$ in $\ell^{\alpha}$ with $\alpha<q$ in each of the following cases:
\begin{itemize}
    \item $n=1$ and $2<q<\infty$.
    \item $n\geq2$ and $1<q<\infty$.
\end{itemize}
This is sharp in the sense that the estimate fails for all $\alpha>q$. 
\end{theorem}
In addition, Bez-Kinoshita-Shiraki \cite{Bez-Kinoshita-Shiraki} also established restricted weak-type estimates for the boundary orthonormal Strichartz estimates at the critical summability exponent $\alpha=q$. It was showed in \cite{Bez-Kinoshita-Shiraki} that, for $n\geq2$, $1<q<\infty$, the estimates
\begin{equation}\label{orthonormal-boundary}
			\bigg\|\sum_{j\in J} \lambda_j|e^{it\Delta}f_j|^2\bigg\|_{L^{q,\infty}(\mathbb{R}, L^\infty(\mathbb{R}^n))}\lesssim \|\{\lambda_j\}_{j\in J}\|_{\ell^{q,1}}, 
\end{equation}	
holds for all families of orthonormal functions $\{f_j\}_{j\in J}$ in $\dot{H}^{s}(\mathbb{R}^n)$ with $s=\frac{n}{2}-\frac{1}{q}$ and all sequences $\{\lambda_j\}_{j\in J}$ in $\ell^{q,1}$. It seems non-trivial to upgrade this to a strong-type estimate.

For $n\geq2$, regarding the summability exponent $\alpha$ for the point $(\frac{1}{p},\frac{1}{q})$ belonging to int $OCDA$, the estimate \eqref{global-orthonormal} holds for $\alpha<q$ and fails for $\alpha>q$. What happens at the critical summability exponent $\alpha=q$ seems to be a delicate matter. As far as we are aware, there are no results available in the literature regarding the summability exponent $\alpha$ for the point $(\frac{1}{p},\frac{1}{q})$ belonging to int $OCDA$. We shall present our new results in this direction below in section \ref{New}.

\subsection{Main new results}\label{New}
Using the frequency-localised estimates in Theorem \ref{frequency-localized} at the critical summability exponent $\alpha=q$, we first establish the following restricted weak-type frequency-global estimates at the same exponent.
\begin{theorem} \label{restricted weak-type}Let $n\geq2$. Let $(\frac{1}{p},\frac{1}{q})$ belonging to int $OCDA$ or $(O,C)$. The estimate 
   \begin{equation*}
			\bigg\|\sum_{j\in J} \lambda_j|e^{it\Delta}f_j|^2\bigg\|_{L^{q,\infty}(\mathbb{R}, L^p(\mathbb{R}^n))}\lesssim \|\{\lambda_j\}_{j\in J}\|_{\ell^{q,1}}, 
\end{equation*}	
holds for all families of orthonormal functions $\{f_j\}_{j\in J}$ in $\dot{H}^{s}(\mathbb{R}^n)$ with $2s=n-\left(\frac{2}{q}+\frac{n}{p}\right)$ and all sequences $\{\lambda_j\}_{j\in J}$ in $\ell^{q,1}$.   
\end{theorem}

We use real interpolation to upgrade the restricted weak-type estimates in Theorem \ref{restricted weak-type} to the strong-type estimates at the critical summability exponent $\alpha=q$.
\begin{theorem}\label{critical strong}
Let $n\geq2$. Let $(\frac{1}{p},\frac{1}{q})$ belonging to int $OCDA$. The estimate 
   \begin{equation*}
			\bigg\|\sum_{j\in J} \lambda_j|e^{it\Delta}f_j|^2\bigg\|_{L^{q}(\mathbb{R}, L^p(\mathbb{R}^n))}\lesssim \|\{\lambda_j\}_{j\in J}\|_{\ell^{q}}, 
\end{equation*}	
holds for all families of orthonormal functions $\{f_j\}_{j\in J}$ in $\dot{H}^{s}(\mathbb{R}^n)$ with $2s=n-\left(\frac{2}{q}+\frac{n}{p}\right)$ and all sequences $\{\lambda_j\}_{j\in J}$ in $\ell^{q}$.   
\end{theorem}
\subsection{Further discussion}
Finally, we make additional comments regarding the remaining interesting issues associated with orthonormal Strichartz estimates for the Schr\"odinger operator.
\begin{itemize}
    \item Regarding the conjectured estimate \eqref{conjecture} at the endpoint $A$ of line segment $[A,D]$, we propose the following analogolus conjecture for the line segment $(O,A]$, as it plays the same role relative to region $OCDA$ as point $A$ does.
    \begin{conjecture}\label{New Conjecture}
       Let $n\geq2$. For any $(\frac{1}{p},\frac{1}{q})$ on the line segment $(O,A]$, the restricted-type estimate
 \begin{equation*}
			\bigg\|\sum_{j\in J} n_j|e^{it\Delta}f_j|^2\bigg\|_{L^q(\mathbb{R}, L^p(\mathbb{R}^n))}\lesssim \|\{n_j\}_{j\in J}\|_{\ell^{q,1}},
\end{equation*}
holds true holds for all families of orthonormal functions $\{f_j\}_{j\in J}$ in $\dot{H}^{s}(\mathbb{R}^n)$ with $2s=n-\left(\frac{2}{q}+\frac{n}{p}\right)$ and all sequences $\{\lambda_j\}_{j\in J}$ in $\ell^{q,1}$.
    \end{conjecture}
In fact, for $n\geq3$, using a semi-classical limiting argument, it was shown by Bez-Hong-Lee-Nakamura-Sawano \cite{Bez-Hong-Lee-Nakamura-Sawano} that if $(\frac{1}{p},\frac{1}{q})$ belongs to $[O,A]$ and $2s=n-\left(\frac{2}{q}+\frac{n}{p}\right)$, then the restricted weak-type estimate
 \begin{equation}\label{interesting}
			\bigg\|\sum_{j\in J} n_j|e^{it\Delta}f_j|^2\bigg\|_{L^{q,\infty}(\mathbb{R}, L^{p,\infty}(\mathbb{R}^n))}\lesssim \|\{n_j\}_{j\in J}\|_{\ell^{q,r}},
\end{equation}
fails for all all families of orthonormal functions $\{f_j\}_{j\in J}$ in $\dot{H}^{s}(\mathbb{R}^n)$ and all sequences $\{\lambda_j\}_{j\in J}$ in $\ell^{q,r}$ for all $r>1$. Whether this statement is true for $n=2$ is a very interesting open problem.
\item For $n=2$, the point $C$ and $D$ coincide. As we know, the single-function Strichartz estimate fails at the Keel-Tao endpoint $D$, and the same for orthonormal Strichartz estimates. Even if the conjecured estimate \eqref{conjecture} were true, deriving orthonormal Strichartz estimates at the critical summability exponent on the line segment $[A,D)$ would remain non-trivial.
\item The possiblity of upgrading the restricted weak-type estimates \eqref{orthonormal-Boundary} and \eqref{orthonormal-boundary} for the boundary case to strong-type versions seems to be an open interesting problem.
\end{itemize}

%%%%%%%%%%%%%%%%%%%%%%%%%%%%%%%%%%%%%%%%%%%%%%%%%%%%%%%%%%%%%%%%%%%%%%%%%
\section{Preliminaries}\label{Preliminaries}
For any appropriate functions $f: \mathbb{R}^n\rightarrow\mathbb{C}$, we use $\hat{f}$ and $\mathcal{F}f$ to denote the Fourier transform of $f$  defined by 
\begin{equation*}
    \hat{f}(\xi)=\int_{\mathbb{R}^n}f(x)e^{-ix\cdot \xi}dx,
\end{equation*}
and the inverse Fourier transform $\mathcal{F}^{-1}f$ and $\check{f}$ defined by
\begin{equation*}
    \check{f}(x)=\frac{1}{(2\pi)^{n}}\int_{\mathbb{R}^n}f(\xi)e^{ix\cdot \xi}d\xi.
\end{equation*}
We define the Schr\"odinger operator $\mathcal{F}(e^{it\Delta}f))(\xi)=e^{-it|\xi|^2}\hat{f}(\xi)$. 

For $p,q\in[1,\infty]$, let $L^p(\mathbb{R}^n)$ denote the usual Lebesgue space and the mixed-norm function space $L^q(\mathbb{R},L^p(\mathbb{R}^n))$ on $\mathbb{R}\times\mathbb{R}^n$ with the norm
\begin{equation*} \|u\|_{L^q(\mathbb{R},L^p(\mathbb{R}^n))}=\left\|\|u(t,\cdot)\|_{L^p(\mathbb{R}^n)}\right\|_{L^q(\mathbb{R})}.
\end{equation*}
Let $\dot{H}^s(\mathbb{R}^n)$ be the homogeneous Sobolev space with the norm
\begin{equation*}
    \|f\|_{\dot{H}^s(\mathbb{R}^n)}=\|(-\Delta)^\frac{s}{2}f\|_{L^2(\mathbb{R}^n)}=\||\cdot|^s\hat{f}\|_{L^2(\mathbb{R}^n)}.
\end{equation*}
\subsection{Littlewood-Paley projections and an interpolation method}\label{LP}
Let $\psi: \mathbb{R}\rightarrow [0,1]$ be a smooth function supported in $[1/2,2]$ such that
\begin{equation*}
    \sum_{k\in\mathbb{Z}}\psi(2^{-k}r)=1,
\end{equation*}
for any $r\neq 0$. For any $k\in\mathbb{Z}$, We define the Littlewood-Paley projection operators $P_k$ on $L^2(\mathbb{R}^n)$ by
\begin{equation*}
    \widehat{P_kf}(\xi)=\psi(2^{-k}|\xi|)\hat{f}(\xi).
\end{equation*}

We introduce the following proposition, which is useful for extending frequency-localized estimates to frequency-global ones. The price we pay is that the resulting frequency-global estimates are in a restricted weak-type form. For further details, we refer the reader to \cite{Bez-Hong-Lee-Nakamura-Sawano} and the key underlying idea traces back to Bourgain \cite{Bourgain}.
\begin{proposition}[\cite{Bez-Hong-Lee-Nakamura-Sawano}]\label{upgrading}
    Let $q_0,q_1>1$, $p,\alpha_0,\alpha_1\geq1$ and $\{g_j\}_{j\in J}$ be a uniformly bounded sequence in $L^{2q_i}(\mathbb{R},L^{2p}(\mathbb{R}^n))$ for each $i=0,1$. If for each $i=0,1,$ there exists $\epsilon_i>0$ such that
    \begin{equation*}
        \left\|\sum_{j\in J} \lambda_j|P_kg_j|^2\right\|_{L^{q_i,\infty}(\mathbb{R},L^{p}(\mathbb{R}^n))}\lesssim 2^{(-1)^{i+1}\epsilon_i k}\|\{\lambda_j\}_{j\in J}\|_{\ell^{\alpha_i}},
    \end{equation*}
    for all $k\in\mathbb{Z}$, then we have
    \begin{equation*}
        \left\|\sum_{j\in J} \lambda_j|g_j|^2\right\|_{L^{q,\infty}(\mathbb{R},L^{p}(\mathbb{R}^n))}\lesssim \|\{\lambda_j\}_{j\in J}\|_{\ell^{\alpha,1}},
    \end{equation*}
    for all sequences $\{\lambda_j\}_{j\in J}\in \ell^{\alpha,1}$, where
    \begin{equation*}
        \frac{1}{q}=\frac{\theta}{q_0}+\frac{1-\theta}{q_1},\;\frac{1}{\alpha}=\frac{\theta}{\alpha_0}+\frac{1-\theta}{\alpha_1} \text{ and } \theta=\frac{\epsilon_1}{\epsilon_1+\epsilon_2}.
    \end{equation*}
\end{proposition}
%%%%%%%%%%%%%%%%%%%%%%%%%%%%%%%%%%%%%%%%%%%%%%%%%%%%%%%%%%%%%%%%%%%%%%%%%
\subsection{Lorentz spaces}
In this subsection, we introduce the Lorentz space $L^{p,r}(\mathbb{R}^n)$ and the Lorentz sequence space $\ell^{p,r}$. We refer the reader to \cite{Stein-Weiss} for further details.

First, let us recall some facts about the rearrangement argument. Let $m$ denote the Lebesgue measure on $\mathbb{R}^n$. Suppose $f$ is a complex-valued measurable function on $\mathbb{R}^n$ such that 
$$m\left(\{x\in\mathbb{R}^n: |f(x)|>t\}\right)=\int_{\{x\in\mathbb{R}^n: |f(x)|>t\}}dx<+\infty,$$ for every $t>0$. The decreasing rearrangement function $f^*$ of $f$ is defined by
\begin{equation*}
f^*(t)=\inf\{s>0: a_f(s)\leq t\},
\end{equation*}
where $a_f$ is the distribution function of $f$ defined by
\begin{equation*}
a_f(t)=m\left(\{x\in\mathbb{R}^n: |f(x)|>t\}\right).
\end{equation*}
Since $f^*$ is decreasing, the maximal function $f^{**}$ of $f^*$ defined by
\begin{equation*}
f^{**}(t)=\frac{1}{t}\int_0^t f^*(s)ds, 
\end{equation*}
is also decreasing and $f^*(t)\leq f^{**}(t)$ for every $t\geq 0$.

For $1\leq p<\infty$ and $1\leq r\leq \infty$, the Lorentz space $L^{p,r}(\mathbb{R}^n)$ is defined as the set of measurable functions $f$ on $\mathbb{R}^n$ such that $\|f\|_{L^{p,r}(\mathbb{R}^n)}<\infty$, where
\begin{equation*}
\|f\|_{L^{p,r}(\mathbb{R}^n)}=
\begin{cases}
\left(\int\limits_0^\infty (t^\frac{1}{p}f^*(t))^r\frac{dt}{t}\right)^\frac{1}{r},&\text{ if } 1\leq r<\infty,\\
\sup\limits_{t>0}t^\frac{1}{p}f^*(t),&\text{ if } r=\infty.
\end{cases}
\end{equation*}
It is well-known that the function $\|\cdot\|_{L^{p,r}(\mathbb{R}^n)}$ defined above is a norm when $r\leq p\neq1$ and a quasi-norm otherwise. In order to obtain a norm in all cases, we define
\begin{equation*}
    \|f\|^*_{L^{p,r}(\mathbb{R}^n)}=
\begin{cases}
\left(\int\limits_0^\infty (t^\frac{1}{p}f^{**}(t))^r\frac{dt}{t}\right)^\frac{1}{r},&\text{ if } 1\leq r<\infty,\\
\sup\limits_{t>0}t^\frac{1}{p}f^{**}(t),&\text{ if } r=\infty.
\end{cases}
\end{equation*}
One can prove that if $1<p<\infty$, $\|\cdot\|^*_{L^{p,r}(\mathbb{R}^n)}$ is a norm, which is equivalent to $\|\cdot\|_{L^{p,r}(\mathbb{R}^n)}$ in the sense that
\begin{equation*}
\|f\|_{L^{p,r}(\mathbb{R}^n)}\leq\|f\|^*_{L^{p,r}(\mathbb{R}^n)}\leq \frac{p}{p-1}\|f\|_{L^{p,r}(\mathbb{R}^n)}, \;\forall f\in L^{p,r}(\mathbb{R}^n).
\end{equation*}
It is clear that $L^{p,p}(\mathbb{R}^n)=L^p(\mathbb{R}^n)$ and the Lorentz spaces are monotone with respect to the second exponent, namely
\begin{equation}\label{embedding}
    L^{p,r_1}(\mathbb{R}^n)\subseteq L^{p,r_2}(\mathbb{R}^n),\;1\leq r_1\leq r_2\leq \infty.
\end{equation}

Finally, we recall the Lorentz sequence space $\ell^{p,r}$. Let $\{\lambda_j\}^\infty_{j=1}\in c_0$ and $\{\lambda^*_j\}^\infty_{j=1}$ be the sequence permuted in a decreasing order. For $1\leq p<\infty$ and $1\leq r\leq \infty$, the Lorentz sequence space $\ell^{p,r}$ is defined as the set of all sequences $\{\lambda_j\}^\infty_{j=1}\in c_0$ such that $\|\{\lambda_j\}^\infty_{j=1}\|_{\ell^{p,r}}<\infty$, where
\begin{equation*}
\|\{\lambda_j\}^\infty_{j=1}\|_{\ell^{p,r}}=
\begin{cases}
\left(\sum\limits_{j=1}^\infty (j^\frac{1}{p}\lambda_j^*)^r\frac{1}{j}\right)^\frac{1}{r},&\text{ if } 1\leq r<\infty,\\
\sup\limits_{j\geq1}j^\frac{1}{p}\lambda_j^*,&\text{ if } r=\infty.
\end{cases}
\end{equation*}
%%%%%%%%%%%%%%%%%%%%%%%%%%%%%%%%%%%%%%%%%%%%%%%%%%%%%%%%%%%%%%%%%%%%%%%
\subsection{Real Interpolation}\label{real-interpolation}
Our proof of upgrading strong-type frequency-global Strichartz estimates involving orthonormal functions from restricted weak-type frequecy-localised versions is based on real interpolation arguments. We shall recall some basic results about the $K$-method of real interpolation. We refer the reader to \cite{Bennett, Bergh-Lofstrom} for details on the development of this theory. Here we only recall the essentials to be used in the sequel. Let $A_0, A_1$ be Banach spaces contained in some larger Banach space $A$. For every $a\in A_0+A_1$ and $t>0$, we define the $K$-functional of real interpolation by
\begin{equation*}
K(t,a,A_0,A_1)=\inf_{a=a_0+a_1}\left(\|a_0\|_{A_0}+\|a_1\|_{A_1}\right).
\end{equation*}For $0<\theta<1$ and $1\leq q\leq \infty$, we denote by $(A_0,A_1)_{\theta,q}$ the real interpolation spaces between $A_0$ and $A_1$ defined as
\begin{equation*}
    (A_0,A_1)_{\theta,q}=\left\{a\in A_0+A_1: \|a\|_{(A_0,A_1)_{\theta,q}}=\left(\int_0^\infty \Big(t^{-\theta}K(t,a,A_0,A_1)\Big)^q\frac{dt}{t}\right)^\frac{1}{q}<\infty\right\}.
\end{equation*}

The real interpolation space identities we shall use are
\begin{equation}\label{L-qp}
    \left(L^{q_0,\infty}(\mathbb{R},L^{p_0}(\mathbb{R}^n)),L^{q_1,\infty}(\mathbb{R},L^{p_1}(\mathbb{R}^n))\right)_{\theta,q}=L^{q}(\mathbb{R},L^{p,q}(\mathbb{R}^n)),
\end{equation}
whenever $\theta\in (0,1)$, $p_0\neq p_1,q_0\neq q_1,p,q\in [1,\infty)$ are such $\frac{1}{q}=\frac{1-\theta}{q_0}+\frac{\theta}{q_1}$ and $\frac{1}{p}=\frac{1-\theta}{p_0}+\frac{\theta}{p_1}$,  and
\begin{equation}\label{l-qp}
    \left(\ell^{q_0,r_0},\ell^{q_1,r_1}\right)_{\theta,r}=\ell^{q,r},
\end{equation}
whenever $\theta\in (0,1)$, $q_0\neq q_1,q\in [1,\infty)$ and $r_0,r_1,r\in[1,\infty]$ are such $\frac{1}{q}=\frac{1-\theta}{q_0}+\frac{\theta}{q_1}$.
%%%%%%%%%%%%%%%%%%%%%%%%%%%%%%%%%%%%%%%%%%%%%%%%%%%%%%%%%%%%%%%%%%%%%%%
\section{Proof of Main new results}
In this section, we give the proof of our new results.
\subsection{Proof of Theorem \ref{restricted weak-type}}
In this section, we give the proof of Theorem \ref{restricted weak-type} by employing Proposition \ref{upgrading} to summ up the frequency-localised estimates in Theorem \ref{frequency-localized} at the critical summability exponent $\alpha=q$.
\begin{proof}
   For the case on the line segment $(O,C)$, it was prove by Bez-Kinoshita-Shiraki \cite{Bez-Kinoshita-Shiraki}. It suffices to prove the case in the interior of $OCDA$. For any fixed point $(\frac{1}{p},\frac{1}{q})\in$ int $OCDA$ and $2s=n-\left(\frac{2}{q}+\frac{n}{p}\right)$, select $\varepsilon>0$ sufficiently small such that $(\frac{1}{p},\frac{1}{q_0})$ and $(\frac{1}{p},\frac{1}{q_1})$ also belong to int $OCDA$, where 
\begin{equation*}
    \frac{1}{q_i}=\frac{1}{q}+(-1)^i\frac{\varepsilon}{2}, \;i=0,1.
\end{equation*}
Then it follows from \eqref{general-k} that 
\begin{align*}
			\bigg\|\sum_{j\in J} \lambda_j|e^{it\Delta}P_k(-\Delta)^{-\frac{s}{2}}f_j|^2\bigg\|_{L^{q_i,\infty}(\mathbb{R},L^{p}(\mathbb{R}^n))}
            &\leq\bigg\|\sum_{j\in J} \lambda_j|e^{it\Delta}P_k(-\Delta)^{-\frac{s}{2}}f_j|^2\bigg\|_{L^{q_i}(\mathbb{R},L^{p}(\mathbb{R}^n))}\\
   &\lesssim 2^{k(n-2s-\frac{2}{q_i}-\frac{n}{p})}\|\{\lambda_j\}_{j\in J}\|_{\ell^{q_i}}\\
   &= 2^{(-1)^{i+1}\varepsilon k}\|\{\lambda_j\}_{j\in J}\|_{\ell^{q_i}},
\end{align*}
holds for any orthonormal system $\{f_j\}_{j\in J} \subseteq L^2(\mathbb{R}^n)$ and $i=0,1$.
Therefore, by Proposition \ref{upgrading}, it immediately implies
\begin{equation*}
			\bigg\|\sum_{j\in J} \lambda_j|e^{it\Delta}(-\Delta)^{-\frac{s}{2}}f_j|^2\bigg\|_{L^{q,\infty}(\mathbb{R}, L^p(\mathbb{R}^n))}\lesssim \|\{\lambda_j\}_{j\in J}\|_{\ell^{q,1}}, 
\end{equation*}	
holds for all families of orthonormal functions $\{f_j\}_{j\in J}$ in $L^2(\mathbb{R}^n)$ and all sequences $\{\lambda_j\}_{j\in J}$ in $\ell^{q,1}$, which is our desired result.
\end{proof}
%%%%%%%%%%%%%%%%%%%%%%%%%%%%%%%%%%%%%%%%%%%%%%%%%%%%%%%%%%%%%%%%%%%%%%%
\subsection{Proof of Theorem \ref{critical strong}}
In this section, we use real interpolation to obtain the optimal strong-type orthonormal Strichartz estimates in the interior $OCDA$, as well as the advantageous condition $q<p$ in this region.
\begin{proof}
 For any fixed point $(\frac{1}{p},\frac{1}{q})$ belonging to int $OCDA$ and $2s=n-\left(\frac{2}{q}+\frac{n}{p}\right)$, we choose $\delta>0$ sufficiently small such that for $(\frac{1}{p_0},\frac{1}{q_0})$ and $(\frac{1}{p_1},\frac{1}{q_1})$ also belong to int $OCDA$, where 
\begin{equation*}
    \frac{1}{p_i}=\frac{1}{p}+(-1)^i\frac{\delta}{n} \quad\text{and}\quad \frac{1}{q_i}=\frac{1}{q}+(-1)^{i+1}\frac{\delta}{2},\quad i=0,1.
\end{equation*}
We can see that $2s=n-\left(\frac{2}{q_i}+\frac{n}{p_i}\right)$ for $i=0,1$, $\frac{1}{p}=\frac{1}{2}\left(\frac{1}{p_0}+\frac{1}{p_1}\right)$ and $\frac{1}{q}=\frac{1}{2}\left(\frac{1}{q_0}+\frac{1}{q_1}\right)$. It follows from Theorem \ref{restricted weak-type} that for $i=0,1$ and all families of orthonormal functions $\{f_j\}_{j\in J}$ in the common space $\dot{H}^{s}(\mathbb{R}^n)$ the estimate
\begin{equation}\label{common}
			\bigg\|\sum_{j\in J} \lambda_j|e^{it\Delta}f_j|^2\bigg\|_{L^{q_i,\infty}(\mathbb{R}, L^{p_i}(\mathbb{R}^n))}\lesssim \|\{\lambda_j\}_{j\in J}\|_{\ell^{q_i,1}}, 
\end{equation}	
holds for all sequences $\{\lambda_j\}_{j\in J}$ in $\ell^{q_i,1}$.  Moreover, by real interpolation space identities \eqref{L-qp} and \eqref{l-qp}, it yields 
  \begin{align*}
       \left(L^{q_0,\infty}(\mathbb{R},L^{p_0,\infty}(\mathbb{R}^n)),L^{q_1}(\mathbb{R},L^{p_1}(\mathbb{R}^n))\right)_{\frac{1}{2},q}&=L^{q}(\mathbb{R},L^{p,q}(\mathbb{R}^n)),\\
        \left(\ell^{q_0,1},\ell^{q_1,1}\right)_{\frac{1}{2},q}&=\ell^{q}.
  \end{align*}
Using real interpolation on the estimates \eqref{common}, we obtain
\begin{equation*}
			\bigg\|\sum_{j\in J} \lambda_j|e^{it\Delta}f_j|^2\bigg\|_{L^{q}(\mathbb{R}, L^{p,q}(\mathbb{R}^n))}\lesssim \|\{\lambda_j\}_{j\in J}\|_{\ell^{q}}.
\end{equation*}	
Note that, for any $(\frac{1}{p},\frac{1}{q})$ belonging to int $OCDA$, it satisfies $q<p$ and by the embedding relation of Lorentz spaces \eqref{embedding}, we have
\begin{equation*}
			\bigg\|\sum_{j\in J} \lambda_j|e^{it\Delta}f_j|^2\bigg\|_{L^{q}(\mathbb{R}, L^{p}(\mathbb{R}^n))}\lesssim \|\{\lambda_j\}_{j\in J}\|_{\ell^{q}},
\end{equation*}	
which are our desired estimates.
\end{proof}
%%%%%%%%%%%%%%%%%%
\section{Strichartz estimates for the kinetic transport equation}
In the final section, we explain the link between orthonormal Strichartz estimates for the Schr\"odinger operator and Strichartz estimates for the velocity average $\rho F$ of the solution $F$ of the kinetic transport equation.  The function $F(x,v,t)=f(x-tv,v)$ satisfies the kinetic transport equation
\begin{equation}
    \begin{cases}
        &(i\partial_t+v\cdot\nabla_x)F(x,v,t)=0, t\in\mathbb{R}, x\in\mathbb{R}^n, v\in\mathbb{R}^n,\\
        &F(x,v,0)=f(x,v),
    \end{cases}
\end{equation}
and 
\begin{equation*}
    \rho F(x,t)=\int_{\mathbb{R}^n} F(x,v,t)dv=\int_{\mathbb{R}^n} f(x-tv,v)dv,
\end{equation*}
is the velocity average of the solution $F$. The kinetic transport equation also enjoys Strichartz estimates of the general form 
\begin{equation}\label{general form}
    \left\|\int_{\mathbb{R}^n} f(x-tv,v)\frac{dv}{|v|^{2s}}\right\|_{L^{q,r}(\mathbb{R}, L^{p,\tilde{r}}(\mathbb{R}^n))}\lesssim \|f\|_{L^{\alpha,\beta}(\mathbb{R}^{2n})},
\end{equation}
for suitable $q, p, \alpha,\beta, r,\tilde{r} \in [1, \infty]$ and $s\in [0,\frac{n}{2})$ satisfying a scaling condition $2s=n-\left(\frac{2}{q}+\frac{n}{p}\right)$.
By a semi-classical limiting argument, the connection between solutions of the free Schr\"odinger equation and the kinetic transport equation is well documented in the following proposition. We refer the reader to Sabin \cite{Sabin} for the special case $s=0$ and Bez-Hong-Lee-Nakamura-Sawano \cite{Bez-Hong-Lee-Nakamura-Sawano} for this general form \eqref{general form}.
\begin{proposition}
Let $q, p\in [1,\infty]$, $\alpha, \beta, r,\tilde{r} \in [1, \infty)$ and $s\in [0,\frac{n}{2})$ be such that $2s=n-\left(\frac{2}{q}+\frac{n}{p}\right)$. If the estimate
\begin{equation*}
			\bigg\|\sum_{j\in J} \lambda_j|e^{it\Delta}f_j|^2\bigg\|_{L^{q,r}(\mathbb{R}, L^{p,\tilde{r}}(\mathbb{R}^n))}\lesssim \|\{\lambda_j\}_{j\in J}\|_{\ell^{\alpha,\beta}}, 
\end{equation*}	
holds for all all families of orthonormal functions $\{f_j\}_{j\in J}$ in $\dot{H}^{s}(\mathbb{R}^n)$ and all sequences $\{\lambda_j\}_{j\in J}$ in $\ell^{\alpha,\beta}$, then the estimate \eqref{general form} holds for any $f\in L^{\alpha,\beta}(\mathbb{R}^{2n})$.
\end{proposition}
As an application of Theorem \ref{globally} and Theorem \ref{critical strong}, we have Strichartz estimates for the solution of the kinetic transport equation with initial data in homogeneous Sobolev spaces.
\begin{theorem}
$(1)$ Let $n\geq 1$. If $(\frac{1}{p},\frac{1}{q})$ belongs to int $OAB$ and $2s=n-\left(\frac{2}{q}+\frac{n}{p}\right)$, then
\begin{equation}\label{General form}
    \left\|\int_{\mathbb{R}^n} f(x-tv,v)\frac{dv}{|v|^{2s}}\right\|_{L^{q}(\mathbb{R}, L^{p}(\mathbb{R}^n)}\lesssim \|f\|_{L^{\alpha^*(p,q)}(\mathbb{R}^{2n})},
\end{equation}
holds for any $f\in L^{\alpha^*(p,q)}(\mathbb{R}^{2n})$.
  
$(2)$ Let $n\geq 2$. If $(\frac{1}{p},\frac{1}{q})$ belongs to int $OCDA$, $2s=n-\left(\frac{2}{q}+\frac{n}{p}\right)$, then the estimate \eqref{General form} holds for any $f\in L^{q}(\mathbb{R}^{2n})$.
\end{theorem}

\section*{Acknowledgments} M. Song would like to express her thanks to Professor Sanghyuk Lee for very enlightening discussions. M. Song is supported by the Guangdong Basic and Applied Basic Research Foundation (Grant No. 2023A1515010656). H. Wu is supported by the National Natural Science Foundation of China (Grant No. 12171399 and 12271041).

\end{document}